\documentclass[12pt]{amsart}
\usepackage{amssymb}
\usepackage{mathabx}
\usepackage{xypic}
\usepackage{hhline}
\usepackage{dcpic}
\usepackage{hyperref}

\newtheorem{theorem}{Theorem}[section]

\newtheorem{hypothesis}[theorem]{Hypothesis}

\newtheorem{property}[theorem]{Property}

\newtheorem{proposition}[theorem]{Proposition}

\newtheorem{remark}[theorem]{Remark}
\newtheorem{definition}[theorem]{Definition}

\numberwithin{theorem}{section}
\numberwithin{equation}{section}

\begin{document}


\makeatletter
\def\Ddots{\mathinner{\mkern1mu\raise\p@
\vbox{\kern7\p@\hbox{.}}\mkern2mu
\raise4\p@\hbox{.}\mkern2mu\raise7\p@\hbox{.}\mkern1mu}}
\makeatother

\newcommand{\OP}[1]{\operatorname{#1}}
\newcommand{\GO}{\OP{GO}}
\newcommand{\leftexp}[2]{{\vphantom{#2}}^{#1}{#2}}
\newcommand{\leftsub}[2]{{\vphantom{#2}}_{#1}{#2}}
\newcommand{\rightexp}[2]{{{#1}}^{#2}}
\newcommand{\rightsub}[2]{{{#1}}_{#2}}
\newcommand{\AI}{\OP{AI}}
\newcommand{\gen}{\OP{gen}}
\newcommand{\prim}{\OP{\star prim}}
\newcommand{\Asai}{\OP{Asai}}
\newcommand{\Image}{\OP{Im}}
\newcommand{\Spec}{\OP{Spec}}
\newcommand{\Ad}{\OP{Ad}}
\newcommand{\Out}{\OP{Out}}
\newcommand{\tr}{\OP{tr}}
\newcommand{\spec}{\OP{spec}}
\newcommand{\scopy}{\OP{end}}
\newcommand{\ord}{\OP{ord}}
\newcommand{\Cent}{\OP{Cent}}
\newcommand{\wellip}{\OP{w-ell}}
\newcommand{\Nrd}{\OP{Nrd}}
\newcommand{\Res}{\OP{Res}}
\newcommand{\temp}{\OP{temp}}
\newcommand{\BC}{\OP{BC}}
\newcommand{\sgn}{\OP{sgn}}
\newcommand{\SU}{\OP{SU}}
\newcommand{\Hom}{\OP{Hom}}
\newcommand{\Inter}{\OP{Int}}
\newcommand{\diag}{\OP{diag}}
\newcommand{\Sym}{\OP{Sym}}
\newcommand{\GSp}{\OP{GSp}}
\newcommand{\GL}{\OP{GL}}
\newcommand{\GSO}{\OP{GSO}}
\newcommand{\bdd}{\OP{bdd}}
\newcommand{\Int}{\OP{Int}}
\newcommand{\art}{\OP{art}}
\newcommand{\vol}{\OP{vol}}
\newcommand{\cusp}{\OP{cusp,\tau}}
\newcommand{\un}{\OP{un}}
\newcommand{\disci}{\OP{disc,\tau_{\it i}}}
\newcommand{\cuspi}{\OP{cusp,\tau_{\it i}}}
\newcommand{\ellip}{\OP{ell}}
\newcommand{\sph}{\OP{sph}}
\newcommand{\gsimp}{\OP{sim-gen}}
\newcommand{\Aut}{\OP{Aut}}
\newcommand{\disc}{\OP{disc,\tau}}
\newcommand{\sdisc}{\OP{s-disc}}
\newcommand{\aut}{\OP{aut}}
\newcommand{\End}{\OP{End}}
\newcommand{\barQ}{\OP{\overline{\mathbf{Q}}}}
\newcommand{\barQp}{\OP{\overline{\mathbf{Q}}_{\it p}}}
\newcommand{\Gal}{\OP{Gal}}
\newcommand{\PGL}{\OP{PGL}}
\newcommand{\simp}{\OP{sim}}
\newcommand{\pri}{\OP{prim}}
\newcommand{\Normal}{\OP{Norm}}
\newcommand{\Ind}{\OP{Ind}}
\newcommand{\St}{\OP{St}}
\newcommand{\unit}{\OP{unit}}
\newcommand{\reg}{\OP{reg}}
\newcommand{\SL}{\OP{SL}}
\newcommand{\Frob}{\OP{Frob}}
\newcommand{\Id}{\OP{Id}}
\newcommand{\GSpin}{\OP{GSpin}}
\newcommand{\Norm}{\OP{Norm}}

\title[Beyond Endoscopic Decomposition]{A weak form of beyond endoscopic decomposition for the stable trace formula of odd orthogonal groups}

\author{Chung Pang Mok}

\address{Department of Mathematics, Purdue University, 150 N. University Street, West Lafayette, IN 47907-2067}

\email{mokc@purdue.edu}

\maketitle

\begin{abstract}
We show that the cuspidal component of the stable trace formula of a special odd orthogonal group over a number field, satisfies a weak form of {\it beyond endoscopic decomposition}. We also study the $r$-stable trace formula, when $r$ is the standard or the second fundamental representation of the dual group. The results are consequences of Arthur's works on endoscopic classification of automorphic representations, together with known results concerning a class of Langlands $L$-functions for special odd orthogonal groups.
\end{abstract}

\section{Introduction}

In the paper [L] (and later expanded in [FLN]), Langlands suggested an approach to establish the Principle of Functoriality for automorphic representations of reductive groups over number fields, an approach that he termed {\it ``beyond endoscopy"}. In this approach, the stable trace formula plays a crucial role; more precisely, one of the main innovation is the insertion of Langlands $L$-functions, or their logarithmic derivatives, into the spectral side of the stable trace formula, via the use of modified form of test functions. An important step in the beyond endoscopy approach is to establish a certain decomposition for the cuspidal component of the stable trace formula, known as beyond endoscopic decomposition.

\bigskip

In this paper, we study the case of a special odd orthogonal group $G =SO(2N+1)$ over a number field $F$, and establish a weak form of the beyond endoscopic decomposition for the cuspidal component of the stable trace formula, by using Arthur's work [A1] on endoscopic classification of automorphic representations for these groups. Equivalently, we would like to put the endoscopic classification framework in [A1] in the setting of the beyond endoscopy framework as in [A2]. In addition, by using results on Langlands $L$-functions $L(s,\pi,r)$, attached to automorphic representation $\pi$ of $G$, with $r$ being either the standard or the second fundamental representation of the dual group $\widehat{G} =Sp(2N,\mathbf{C})$, we show that the $r$-stable trace formula for $G$ (in the sense of [A2]) exists, and that it also satisfies a similar weak form of beyond endoscopic decomposition; in the case when $r$ is the standard representation, the $r$-stable trace formula in fact vanishes identically.

\bigskip

To state these results, we must recall some formulation of stable trace formula of Arthur and his work on endoscopic classification, which we will turn to in the next section. The weak beyond endoscopic decomposition of the cuspidal component of the stable trace formula will be established in section 3. In section 4 we study the $r$-stable trace formula for $r$ being the standard of the second fundamental representations of $\widehat{G}$ The main results are established as Theorem 3.6, Theorem 4.1 and Theorem 4.4. In the final section 5 we will pose some closely related questions.

\section*{Acknowledgement}

The author would like to thank Professor Jiang for informing him in 2013 of his paper [J], which formed the inspiration for the present paper. He would also like to thank Professor Shahidi for various encouragements and discussions during the publication process of the paper. The author is grateful to the referee for the careful reading of the paper.

\section{R\'esum\'e on stable trace formula and endoscopic classification}

In the general context of the Arthur trace formula, one has a connected reductive group $G$ over a number field $F$. In the present paper, we will need the case where $G$ is the split form of the special odd orthogonal group $SO(2N+1)$ over $F$ (with $N \geq 1$), or finite products of such groups. The Hecke space of test functions on the adelic group $G(\mathbf{A}_F)$ will be denoted as 
\[
\mathcal{H}(G(\mathbf{A}_F)) = \rightexp{\prod_v}{\prime} \mathcal{H}(G(F_v))
\]
\noindent which is a restricted direct product of the local Hecke space $\mathcal{H}(G(F_v))$. Here the restricted direct product is taken with respect to the unit element of the spherical Hecke algebra $\mathcal{H}^{\sph}(G(F_v))$ at places $v$ of $F$. Then for each non-negative real number $\tau$, one has the linear form $I^G_{\disc}$ on $\mathcal{H}(G(\mathbf{A}_F))$, known as the discrete component of the invariant trace formula for $G$; here and thereafter, the parameter $\tau$ controls the norm of the
imaginary part of the archimedean infinitesimal characters of representations. We refer the reader to section 3.1 of [A1] for more detailed discussions. The value of the linear form $I^G_{\disc}$ on a test function $f \in \mathcal{H}(G(\mathbf{A}_F))$ depends only on the invariant orbital integral $f_G$, with $f_G$ being regarded as a function on the set of regular semi-simple conjugacy classes of $G(\mathbf{A}_F)$.

\bigskip

For questions related to the Principle of Functoriality, it is the stable trace formula that plays the crucial role [A4,A5,A6]. More precisely, we have the linear form $S^G_{\disc}$ on $\mathcal{H}(G(\mathbf{A}_F))$, known as the discrete component of the stable trace formula for $G$. The linear form satisfies the important condition of being stable, namely that its value on a test function $f \in \mathcal{H}(G(\mathbf{A}_F))$ depends only on the stable orbital integral $f^G$, with $f^G$ being regarded as a function on the set of regular semi-simple stable conjugacy classes of $G(\mathbf{A}_F)$; for these notions we refer for example to section 2.1 and section 3.2 of [A1]. We denote by $\mathcal{S}(G(\mathbf{A}_F))$ the space spanned by $f^G$ for $f \in \mathcal{H}(G(\mathbf{A}_F))$. One has a corresponding local definition, and we have
\[
\mathcal{S}(G(\mathbf{A}_F)) = \rightexp{\bigotimes_v}{\prime} \mathcal{S}(G(F_v))
 \]
\noindent here the restricted tensor products is taken with respect to the stable orbital integral of the unit element of the spherical Hecke algebra at places $v$ of $F$. 
\bigskip

Following [A1], given a linear form $S$ on $\mathcal{H}(G(\mathbf{A}_F))$ that is stable, we will denote by $\widehat{S}$ the corresponding linear form on the space of stable orbital integrals $\mathcal{S}(G(\mathbf{A}_F))$; in other words, for $f \in \mathcal{H}(G(\mathbf{A}_F))$,
\[
S(f) = \widehat{S}(f^G)
\]

\bigskip

The stable linear form $S^G_{\disc}$ is constructed from $I^G_{\disc}$ via the inductive procedure as described for example in section 3.2 of [A1]. This inductive procedure critically relies on the existence of the Langlands-Shelstad transfer and the fundamental lemma. As is well known, these are now all established theorems; for the precise reference we refer to section 2.1 of [A1]. Finally, the linear forms $I^G_{\disc}, S^G_{\disc}$ satisfy an admissible condition, {\it c.f.} section 3.1 of [A1].

\bigskip

We now specialize to the case where $G$ is the split $SO(2N+1)$. In the work [A1], Arthur gave explicit formula for $S^G_{\disc}$, in terms of the set of (formal) Arthur parameters $\Psi(G)$ for the group $G$. More precisely, since we have fixed the parameter $\tau$, it suffices to work with the subset
\[
\Psi^{\tau}(G) \subset \Psi(G)
\]
consisting of parameters $\psi$, having the property that, the localizations $\{ \psi_v \}_{v | \infty}$ of $\psi$ at the archimedean places parametrizes archimedean $L$-packets consisting of representations of $\prod_{v | \infty} G(F_v)$ with the norm of the
imaginary part of their infinitesimal characters being equal to $\tau$. Then for each such parameter $\psi \in \Psi^{\tau}(G)$, one has a corresponding stable linear form $S^G_{\disc,\psi}$, and one has a decomposition:
\[
S^G_{\disc} = \sum_{\psi \in \Psi^{\tau}(G)} S^G_{\disc,\psi}.
\]  
Since the linear form $S^G_{\disc}$ is admissible, ({\it c.f.} section 3.1 of [A1]), one has that, for $f \in \mathcal{H}(G(\mathbf{A}_F))$, the sum:
\[
S^G_{\disc} (f)= \sum_{\psi \in \Psi^{\tau}(G)} S^G_{\disc,\psi}(f).
\]
is a finite sum, and the set of indices $\psi$ that contribute non-trivially to the sum, depends only on a choice of Hecke type for $f$ ({\it c.f.} {\it loc. cit.}).

\bigskip

In [A1], Arthur showed that each $\psi \in \Psi^{\tau}(G)$ parametrizes a set of unitary representations of $G(\mathbf{A}_F)$, known as the global Arthur packet attached to $\psi$. Furthermore, the $L^2$- discrete automorphic spectrum of $G$:
\[
L^2_{\disc}(G(F) \backslash G(\mathbf{A}_F))
\]
has a decomposition in terms of the global Arthur packets parametrized by the subset $\Psi^{\tau}_{2}(G) \subset \Psi^{\tau}(G)$ consisting of square-integrable parameters. We will be particularly concerned with the subset $\Phi^{\tau}_{2}(G) \subset \Psi^{\tau}_{2}(G)$, consisting of parameters that are both square-integrable and generic. A parameter $\phi$ that belongs to $\Phi^{\tau}_{2}(G)$ will be referred to as cuspidal. 

\bigskip

We refer the reader to chapter one of [A1] for these notions. For our purpose, it suffices to recall that in [A1], a cuspidal parameter $\phi \in \Phi^{\tau}_{2}(G)$ is the data given as an unorderd formal direct sum:
\begin{eqnarray}
\phi = \phi_1 \boxplus \cdots \boxplus \phi_k
\end{eqnarray}
with each $\phi_i$ ($i=1,\cdots,k$) being an unitary cuspidal automorphic representation on $GL(m_i,\mathbf{A}_F)$ that is of symplectic type; here each $m_i$ is an even integer, and one requires that $2N =\sum_{i=1}^k m_i$, and that the $\phi_i$'s are mutually non-equivalent; here for $1 \leq s,t \leq k$, we say that $\phi_s$ and $\phi_t$ are equivalent is $m_s=m_t$, and $\phi_s,\phi_t$ are equivalent as irreducible representations of $\GL(m_s,\mathbf{A}_F)$. Finally $\phi_i$ being of symplectic type means that the exterior square $L$-function $L(s,\phi_i,\Lambda^2)$ has a pole at $s=1$ (in which case the pole is a simple pole). This definition is taken to be the formal substitute of the notion of global Arthur parameters, in the absence of the global automorphic Langlands group; {\it c.f.} chapter 1 of [A1] for more information. Finally, we will refer to a cuspidal parameter $\phi \in \Phi^{\tau}_{2}(G)$ as simple, if $k=1$ in (2.1) above. We will denote the subset of $\Phi_2^{\tau}(G)$ consisting of simple cuspidal parameters of $G$ as $\Phi^{\tau}_{\simp}(G)$.

\bigskip

We now define the following stable linear form:
\begin{eqnarray}
S^G_{\cusp} :=\sum_{\phi \in \Phi^{\tau}_2(G)} S^G_{\disc,\phi} 
\end{eqnarray}
which we will refer to as the cuspidal component of the stable trace formula for $G$. It plays the main role in the formulation of the results in section 3 and 4.

\begin{remark}

\end{remark}

\noindent A consequence of the work [A1] is that, on assuming the Generalized Ramanujan Conjecture for general linear groups, one has the validity of the corresponding Generalized Ramanujan Conjecture for $G$: a discrete automorphic representation $\pi$ of $G(\mathbf{A}_F)$ is tempered at every place of $F$, if and only if it belongs to a global Arthur packet parameterized by a cuspidal parameter. Thus we see that our definition of the linear form $S^G_{\cusp}$ is consistent with that given in section 2 of [A2].

\bigskip

We similarly define the following stable linear form:
\begin{eqnarray}
\leftexp{\star}{P}^G_{\cusp} :=\sum_{\phi \in \Phi^{\tau}_{\simp}(G)} S^G_{\disc,\phi}
\end{eqnarray}
which we will refer to as the $\star$-primitive component of the linear form $S^G_{\disc}$. 

\begin{remark}
\end{remark}

\noindent The linear forms $S^G_{\cusp}$ and $\leftexp{\star}{P}^G_{\cusp}$ are again admissible ({\it c.f.} section 3.1 of [A1]). Note that their definition invokes the formalism of the endoscopic classification. The definition of $\leftexp{\star}{P}^G_{\cusp}$ is intended to be the analogue of the conjectural primitive component $P^G_{\cusp}$ of $S^G_{\disc}$, {\it c.f.} section 2 of [A2], and also the discussions in item VI of section 5 below. In fact, the construction of the primitive component $P^G_{\cusp}$ of $S^G_{\disc}$, whose existence is still conjectural at this moment, would already be an important step in Langlands' beyond endoscopy apporach. 

\bigskip

One of the main results in [A1] is an explicit formula for the stable linear form $S^G_{\disc,\psi}$, for any parameter $\psi \in \Psi^{\tau}(G)$; this result is known as the stable multiplicity formula. For our purpose we only need the stable multiplicity formula for the cuspidal parameters $\Phi^{\tau}_{2}(G)$:

\begin{proposition} ({\it c.f.} Theorem 4.1.2 of [A1])
For a parameter $\phi \in \Phi^{\tau}_{2}(G)$, with 
\[
\phi = \phi_1 \boxplus \cdots \boxplus \phi_k
\]
as in (2.1) above, and $f \in \mathcal{H}(G(\mathbf{A}_F))$, we have the formula
\begin{eqnarray}
S^G_{\disc,\phi}(f) = \frac{1}{|\mathcal{S}_{\phi}|} \cdot f^G(\phi). 
\end{eqnarray}
Here $|\mathcal{S}_{\phi}|=2^{k-1}$, and the term $f^G(\phi)$ is value of the stable character, evaluated at the test function $f$, attached to the global packet parametrized by $\phi$.
\end{proposition}

\begin{remark}
\end{remark}

\noindent $|\mathcal{S}_{\phi}|$ is the order of the group $\mathcal{S}_{\phi}$, defined to be 
\begin{eqnarray*}
& & \mathcal{S}_{\phi} :=\pi_0(\overline{S}_{\phi}) \\
& & \overline{S}_{\phi} = S_{\phi} /Z(\widehat{G})
\end{eqnarray*}
and $S_{\phi}$ being a reductive (non-connected in general) subgroup of $\widehat{G}$, as defined for example in chapter 1 of [A1].
\bigskip

In particular we have the following formula, for $f \in \mathcal{H}(G(\mathbf{A}_F))$:
\begin{eqnarray}
\leftexp{\star}{P}^G_{\cusp}(f) = \sum_{\phi \in \Phi^{\tau}_{\simp}(G)} f^G(\phi)
\end{eqnarray}

\section{Weak Beyond Endoscopic Decomposition}

We now establish the weak form of the beyond endoscopic decomposition, which we will refer to as $\star$-beyond endoscopic decomposition of $S^G_{\cusp}$, in this section. First we make the definition:
\begin{definition}
An elliptic $\star$-beyond endscopic data of $G$ is a pair $(H,\rho)$, where $H$ is isomorphic to a product of special odd orthogonal groups over $F$:
\begin{eqnarray}
& & H \cong H_1 \times \cdots \times H_k \\
& & H_i = SO(m_i +1) \nonumber
\end{eqnarray}
with the $m_i$ being positive even integers, satisfying $2N = \sum_{i=1}^k m_i$, and $\rho$ being the datum of an embedding:
\[
\widehat{H} \hookrightarrow \widehat{G}
\] 
\end{definition}
Here and in what follows, we always work with the standard based root data of the orthogonal groups involved. Two such pairs $(H_1,\rho_1)$ and $(H_2,\rho_2)$ are said to be equivalent, if there exists an isomorphism $\alpha: H_1 \cong H_2$ over $F$, and $g \in \widehat{G}$, such that for the dual isomorphism $\widehat{\alpha}: \widehat{H}_1 \cong \widehat{H}_2$, we have $\rho_2 \circ \widehat{\alpha}$ and $\rho_1$ differ by conjugation by $g$. It can be seen that the equivalence class of the pair $(H,\rho)$ as above is completely characterized by the unordered partition $2N =\sum_{i=1}^k m_i$, so is characterized by the underlying group $H$. We will often abuse the notation and denote an equivalence class of elliptic $\star$-beyond endoscopic data simply as $(H,\rho)$ or even as $H$.

\bigskip

\noindent Here the term ``elliptic" refers to the condition that the cardinality of the group
\begin{eqnarray}
\overline{S}_{\rho} := S_{\rho} /Z(\widehat{G}) 
\end{eqnarray}
is finite, where
\begin{eqnarray}
S_{\rho} :=\Cent_{\widehat{G}}(Im \rho)
\end{eqnarray}
In fact the cardinality $|\overline{S}_{\rho}|$ equal to the number $2^{k-1}$ for a pair $(H,\rho)$ as in (3.1) above. We put
\begin{eqnarray}
\iota(G,H) := \big|\overline{S}_{\rho} \big|^{-1}
\end{eqnarray}

\bigskip
Although we do not need this, one can define a $\star$-beyond endoscopic data for $G$ to be given by a pair $(M,\rho_M)$, where $M$ is a Levi subgroup of $H$ for some elliptic $\star$-beyond endoscopic data $(H,\rho)$ as above, and $\rho_M$ being the embedding $\widehat{M} \hookrightarrow \widehat{G}$ given by the composition of the dual embedding $\widehat{M} \hookrightarrow \widehat{H}$ with the embdding $\rho$ (and with a similar definition of equivalence of such pairs). One defines the groups $S_{\rho_M}, \overline{S}_{\rho_M}$ exactly as in (3.2), (3.3) above. Then one sees that a $\star$-beyond endoscopic data $(M,\rho_M)$ is elliptic exactly when $|\overline{S}_{\rho_M}|$ is finite.

\bigskip

\begin{remark}
\end{remark}

\noindent In the setting of standard endoscopy, an elliptic endoscopic data for $G$ consists of pairs $(H,\rho)$ as in (3.1) above, but only with $k=1,2$, so it is more restrictive.

\begin{remark}
\end{remark}
Again, definition 3.1 is only a very special case of the general definition of the notion of (elliptic) beyond endoscopic data, {\it c.f.} section 2 of [A2].

\bigskip

For an elliptic $\star$-beyond endoscopic data $(H,\rho)$ of $G$ as above, we may assume without loss of generality that
\[
H= H_1 \times \cdots \times H_k
\]
with $H_i = SO(m_i+1)$ for positive even integers $m_i$ as before. Put 
\[
\Phi^{\tau}_2(H) := \bigcup_{\substack{\tau_1,\cdots,\tau_k \\ \Sigma_{i=1}^k \tau_i^2=\tau^2 } } \prod_{i=1}^k \Phi^{\tau_i}_2(H_i). 
\]
Note that if $\phi \in \Phi^{\tau}_2(G) - \Phi^{\tau}_{\simp}(G)$, then $\phi$ arises from a $\phi_H \in \Phi^{\tau}_2(H)$, for an elliptic $\star$-beyond endoscopic data $(H,\rho)$ of size smaller than $G$. In the framework of [A1] we can say that $\phi$ is the composition of $\phi_H$ with $\rho$, and hence $\phi$ should not be regarded as $\star$-primitive with respect to $G$, which gives the motivation for the definition of $\leftexp{\star}{P}^G_{\cusp}$.

\bigskip

We now define $\Phi^{\tau}_{\prim}(H)$, the set of $\star$-primitive parameters of $H$, to be the subset of $ \bigcup_{\substack{\tau_1,\cdots,\tau_k \\ \Sigma_{i=1}^k \tau_i^2=\tau^2 }} \prod_{i=1}^k \Phi^{\tau_i}_{\simp}(H_i) \subset \Phi^{\tau}_2(H)$, consisting of elements:
\begin{eqnarray}
\phi_H = \phi_1 \times \cdots \times \phi_k
\end{eqnarray}
with each $\phi_i \in \Phi^{\tau_i}_{\simp}(H_i) =\Phi^{\tau_i}_{\simp}(SO(m_i+1))$, mutually inequivalent. In particular $\Phi^{\tau}_{\prim}(H_i) =\Phi^{\tau}_{\simp}(H_i)$ for each factor $H_i$.

\bigskip
Here we are implicitly using the fact that the norm of the imaginary part of archimedean infinitesimal characters is preserved under Langlands functorial transfer at the archimedean places, and so one has compatibility of the parameter $\tau$ between $G$ and the various $(H,\rho)$.

\begin{remark}
\end{remark}
\noindent A motivation for this definition is as follows: indeed if one considers a parameter $\phi_H = \phi_1 \times \cdots \times \phi_k \in  \bigcup_{\substack{\tau_1,\cdots,\tau_k \\ \Sigma_{i=1}^k \tau_i^2=\tau^2 } } \prod_{i=1}^k \Phi^{\tau_i}_{\simp}(H_i) \subset \Phi^{\tau}_2(H)$, say with $\phi_s \cong \phi_t$ for some $ 1 \leq s < t \leq k$ (in particular $m_s=m_t$ and $H_s=H_t$), then we can regard $\phi_H$ as arising from the smaller group:
\[
H^t := H_1 \times \cdots \times H_{t-1} \times H_{t+1} \times \cdots \times H_k
\]
with respect to the partial diagonal embedding
\begin{eqnarray*}
\widehat{H^t} = \widehat{H}_1 \times \cdots \times \widehat{H}_{t-1} \times \widehat{H}_{t+1}  \times  \cdots  \times   \widehat{H}_k  \hookrightarrow  \widehat{H} =  \widehat{H}_1 \times \cdots \times \widehat{H}_k 
\end{eqnarray*}
\begin{eqnarray*}
(g_1,\cdots , g_{t-1},g_{t+1},\cdots, g_k) \mapsto  (g_1,\cdots , g_{t-1},g_s,g_{t+1},\cdots, g_k)
\end{eqnarray*}
and hence such a $\phi_H$ should not be regarded as $\star$-primitive with respect to $H$.
\bigskip

We then define the stable linear form on $\mathcal{H}(H(\mathbf{A}_F))$:

\begin{eqnarray}
\leftexp{\star}{P}^H_{\cusp}(f^{\prime}) = \sum_{\phi_H \in \Phi^{\tau}_{\prim}(H)} (f^{\prime})^H(\phi_H), \,\ f^{\prime} \in \mathcal{H}(H(\mathbf{A}_F))
\end{eqnarray}
And just as $\leftexp{\star}{P}^G_{\cusp}$, the linear form $\leftexp{\star}{P}^H_{\cusp}$ is admissible.

\begin{remark}
\end{remark}

\noindent For a test function $f^{\prime}  \in \mathcal{H}(H(\mathbf{A}_F)) = \prod_{i=1}^r \mathcal{H}(H_i(\mathbf{A}_F))$, with $f_i^{\prime}$ being the component of $f$ in the $i$-th factor, we have:
\begin{eqnarray*}
S^H_{\disc}(f^{\prime}) = \sum_{\substack{\tau_1,\cdots,\tau_k  \\ \Sigma_{i=1}^k \tau_i^2=\tau^2 }}  \big( \prod_{i=1}^r S^{H_i}_{\disci}(f^{\prime}_i) \big)
\end{eqnarray*}
but in general
\[
\leftexp{\star}{P}_{\cusp}^H(f^{\prime}) \neq  \sum_{\substack{ \tau_1,\cdots,\tau_k  \\ \Sigma_{i=1}^k \tau_i^2=\tau^2 } } \big( \prod_{i=1}^r \leftexp{\star}{P}_{\cuspi}^{H_i}(f^{\prime}_i) \big).
\]
precisely by the reason explained as above, namely that $\Phi^{\tau}_{\prim}(H)$ is in general a proper subset of $ \bigcup_{\substack{ \tau_1,\cdots,\tau_k  \\  \Sigma_{i=1}^k \tau_i^2=\tau^2 }} \prod_{i=1}^k \Phi^{\tau_i }_{\simp}(H_i)$.

\bigskip
We can now state the main result of this section (here below recall that the factor $\iota(G,H)$ is defined in (3.4)):
\begin{theorem}
We have the $\star$-beyond endoscopic decomposition:
\begin{eqnarray}
S^G_{\cusp}(f) = \sum_{(H,\rho)} \iota(G,H) \cdot \widehat{\leftexp{\star}{P}^H_{\cusp}}(f^H), \,\ f \in \mathcal{H}(G(\mathbf{A}_F))
\end{eqnarray}
Here $(H,\rho)$ runs over the equivalence classes of elliptic $\star$-beyond endoscopic data of $G$. Finally $f^H$ is the stable transfer of $f$ to $H$, which we will define below. 
\end{theorem}

\bigskip

Here given an elliptic $\star$-beyond endoscopic data $(H,\rho)$ of $G$, the stable transfer is a map:
\begin{eqnarray}
f \rightarrow f^H, \,\ f \in \mathcal{H}(G(\mathbf{A}_F)).
\end{eqnarray}
the image $f^H$ is an element of $\mathcal{S}(H(\mathbf{A}_F))$, and it depends only on the image $f^G$ of $f$ in $\mathcal{S}(G(\mathbf{A}_F))$. It is defined purely locally, namely that given 
\[
f=\rightexp{\prod_v}{\prime} f_v
\]
one has 
\[
f^H = \rightexp{\bigotimes_v}{\prime} f_v^H
\]
So it suffices to work locally. By the results on local Langlands correspondence for $G$ as established in chapter 6 of [A1], the set of bounded local Langlands parameters for $G$ (resp. $H$) over $F_v$ parametrizes tempered $L$-packets of $G(F_v)$ (resp. $H(F_v)$), and the space of stable orbital integrals on regular semi-simple elements of $G(F_v)$ (resp. $H(F_v)$) corresponds to the Paley-Wiener space on the set of bounded local Langlands parameters for $G$ (resp. $H$) over $F_v$, under the map given by taking stable characters.

\bigskip
As before, write $H = H_1 \times \cdots \times H_k$, with $H_i = SO(m_i +1)$. To construct the stable transfer, let $f_v \in \mathcal{H}(G(F_v))$. Firstly consider the following function $J$ on the set of bounded local Langlands parameters of $G$ over $F_v$, determined by the condition: for each bounded local Langlands parameter 
\[
\phi_{v}: L_{F_v} \rightarrow \widehat{G}
\]
of $G$ over $F_v$ (here $L_{F_v}$ being the local Langlands group of $F_v$), we have:
\begin{eqnarray}
J(\phi_{v}) = f_v^G(\phi_{v}).
\end{eqnarray}
The function $J$ belongs to the Paley-Wiener space on the set of bounded local Langlands parameters of $G$ over $F_v$.  Now define the function $J^H$ on the set of bounded local Langlands parameters of $H$ over $F_v$, determined by the condition: for each bounded local Langlands parameter 
\[
\phi_{H,v}: L_{F_v} \rightarrow \widehat{H}
\]
of $H$ over $F_v$, we have:
\begin{eqnarray*}
J^H(\phi_{H,v}) = J( \rho \circ \phi_{H,v}).
\end{eqnarray*}

From the fact that $\rho: \widehat{H} \hookrightarrow \widehat{G}$ is an embedding that induces isomorphism on the corresponding maximal tori, we see that the function $J^H$ belongs to the Paley-Wiener space on the set of bounded local Langlands parameters of $H$ over $F_v$. Thus there exists unique $f_v \in \mathcal{S}(H(F_v))$ such that
\[
f_v^H(\phi_{H,v}) =J(\phi_{H,v})
\]
and at the same time, we see that the element $f_v^H \in \mathcal{S}(H(F_v))$ is uniquely determined by $f_v$ (and in fact, depends only on the image $f_v^G$ of $f_v$ in $\mathcal{S}(G(F_v))$). The map $f_v \mapsto f_v^H$ is the stable transfer that we look for.

\bigskip

A result of independent interest is that, in our present setting, the stable transfer as described above can also be constructed in terms of the Kottwitz-Shelstad transfer in the setting of twisted endoscopy. We refer to section 6 for this result, which will not be needed in the main part of the paper. 

\bigskip

\bigskip

In addition, from the defining condition (3.9) in terms of local Langlands parameters, one sees that, if $f_v$ belongs to the spherical Hecke algebra $\mathcal{H}^{\sph}(G(F_v))$, then $f_v^H$ is given by $b(f_v)^H$; here $b(f_v) \in \mathcal{H}^{\sph}(H(F_v))$ is the image of $f_v$ under the algebra homomorphism:
\begin{eqnarray}
b: \mathcal{H}^{\sph}(G(F_v)) \rightarrow \mathcal{H}^{\sph}(H(F_v))
\end{eqnarray}
that is induced, via the Satake isomorphism, by the embedding $\rho: \widehat{H} \hookrightarrow \widehat{G}$; in particular $b$ takes the unit element of $\mathcal{H}^{\sph}(G(F_v))$ to the unit element of $\mathcal{H}^{\sph}(H(F_v))$. Thus we see that the stable transfer map
\[
\mathcal{H}(G(\mathbf{A}_F)) \rightarrow \mathcal{S}(H(\mathbf{A}_F))
\] 
is defined, and factors through $\mathcal{S}(G(\mathbf{A}_F))$.

\bigskip

Now back to the global context. Given $\phi \in \Phi^{\tau}_2(G)$, there is, up to equivalence, an unique elliptic $\star$-beyond endoscopic data $(H,\rho)$ of $G$, such that $\phi$ is the composition of $\phi_H$ with $\rho$, for some $\phi_H \in \Phi^{\tau}_{\prim}(H)$, satisfying $|\mathcal{S}_{\phi}| =|\overline{S}_{\rho}|$; as a matter of fact, in the context of chapter one of [A1], the definition of the group $S_{\phi}$ coincides exactly with $S_{\rho}$, hence $\mathcal{S}_{\phi} = \pi_0(\overline{S}_{\phi})=\pi_0(\overline{S}_{\rho}) =\overline{S}_{\rho}$ (as $\overline{S}_{\rho}$ is already finite). Now given
\[
f = \rightexp{\prod_v}{\prime} f_v \in \mathcal{H}(G(\mathbf{A}_F)),
\]
we have:
\begin{eqnarray*}
& & f^H(\phi_H) \\&=& \prod_v f_v^H(\phi_{H,v})
= \prod_v f_v^G(\rho \circ \phi_{H,v}) \\ &=& \prod_v f_v^G(\phi_v)=  f^G(\phi)
\end{eqnarray*}

The proof of theorem 3.6 now follows, by using the stable multiplicity formula (Proposition 2.3), from the following computation:
\begin{eqnarray*}
& & S^G_{\cusp}(f) = \sum_{\phi \in \Phi^{\tau}_2(G)} S^G_{\disc,\phi}(f) \\
&=& \sum_{\phi \in \Phi^{\tau}_2(G)}  \frac{1}{|\mathcal{S}_{\phi}|} \cdot  f^G(\phi) \\
&=& \sum_{(H,\rho)} \frac{1}{|\overline{S}_{\rho} |} \sum_{\phi_H \in \Phi^{\tau}_{\prim}(H)} f^H(\phi_H) \\&=&  \sum_{(H,\rho)} \iota(G,H) \cdot \widehat{\leftexp{\star}{P}^H_{\cusp}}(f^H).
\end{eqnarray*}

\section{The $r$-stable trace formula}

We now give examples of $r$-stable trace formula (or more simply $r$-trace formula), in the sense of [A2]. As in the introduction, $r$ will be a finite dimensional algebraic representation of $\widehat{G}=Sp(2N,\mathbf{C})$. We will be particularly concerned with the case where $r$ is a fundamental representation of $\widehat{G}$; more precisely, denote by $std$ the standard representation of $\widehat{G}$ (of dimension equal to $2N$). Then for $a=1,\cdots,N$, the $a$-th fundamental representation $r_a$ on $\widehat{G}$ fits into the following split short exact sequence:
\[
0 \rightarrow r_a \rightarrow \Lambda^a std \rightarrow \Lambda^{a-2} std \rightarrow 0
\]
{\it c.f.} the discussion in section 2 of [J]. In particular $r_1 = std$, and $r_2$ is of dimension $2N^2-N-1$ with $\Lambda^2 std = r_2 \oplus 1$.

\bigskip

In Langlands' beyond endoscopy proposal, one would like to construct a certain limiting form of the stable trace formula for $G$, known as the $r$-stable trace formula, for each algebraic representation $r$ of $\widehat{G}$. We briefly recall the situation, {\it c.f.} section 2 [A2] and the afterword section of [A3]. Given a test function $f \in \mathcal{H}(G(\mathbf{A}_F))$, fix a finite set $S$ of places of $F$ (we assume that $S$ contains all the archimedean places of $F$) such that the component of $f$ at places $v \notin S$ belong to the spherical Hecke algebra $\mathcal{H}^{\sph}(G(F_v))$. Then for each place $w \notin S$, and $n \geq 1$, denote by $\leftexp{n}{h}^r_w \in \mathcal{H}^{\sph}(G(F_w))$ the function, whose Satake transform $\widehat{  \leftexp{n}{h}^r_w}$ satisfies: 
\begin{eqnarray}
\widehat{\leftexp{n}{h}^r_w(c) }= tr (r(c^n))
\end{eqnarray}
for every semi-simple conjugacy class $c$ of $\widehat{G}$. For $n=1$ we simply denote the function as $h^{r}_w$. We then replace $f$ by the function $\leftexp{n}f^{r,w} \in \mathcal{H}(G(\mathbf{A}_F))$, where the local component of $\leftexp{n}{f}^{r,w}$ at each place $v \neq w$ is the same as $f$, while at the place $w$, the local component of $\leftexp{n}{f}^{r,w}$ is equal to the product (in $\mathcal{H}^{\sph}(G(F_w))$) of $\leftexp{n}{h}^r_w$ with the local component $f_w$ of $f$ at $w$. Again when $n=1$ we simply denote this as $f^{r,w}$.

\bigskip

Langlands [L] proposed to establish the existence of the limit:
\begin{eqnarray}
 S^{r,G}_{\cusp}(f):=\lim_{X \rightarrow \infty }\frac{1}{X} \sum_{Nw \leq X, w \notin S}  \log Nw \cdot S^{G}_{\cusp}(f^{r,w})
\end{eqnarray}
One sees that the limit does not depend on the choice of $S$ if it exists. If the limit exists, then we refer to the linear form $S^{r,G}_{\cusp}$ on $\mathcal{H}(G(\mathbf{A}_F))$ as the $r$-stable trace formula.

\bigskip
One can also define a variant of (4.2), as the residue at $s=1$ of a Dirichlet series $\mathcal{S}^{r,G}_{\cusp}(f,s)$:
\begin{eqnarray}
& & \mathfrak{S}^{r,G}_{\cusp}(f) := \lim_{s \rightarrow 1} (s-1) \cdot \mathcal{S}^{r,G}_{\cusp}(f,s) \\
& & \mathcal{S}^{r,G}_{\cusp}(f,s) := \sum_{ w \notin S}  \sum_{n \geq 1} \frac{\log Nw}{Nw^{ns}} \cdot S^{G}_{\cusp}(\leftexp{n}{f}^{r,w}) \nonumber
\end{eqnarray}
with $s $ a complex variable. Here we interpret the limit in (4.3) as requiring that the Dirichlet series $\mathcal{S}^{r,G}_{\cusp}(f,s)$ converges absolutely in some right half plane, and possesses analytic continuation holomorphic in the region $Re(s) >1$. If the limit (4.3) exists, we will denote it as $\mathfrak{S}^{r,G}_{\cusp}(f)$, and we also refer to it as the $r$-stable trace formula (it is expected that when the limits $S^{r,G}_{\cusp}(f),\mathfrak{S}^{r,G}_{\cusp}$ both exist, they are equal). Finally, the script font used for $\mathcal{S}^{r,G}_{\cusp}$ is to suggest that it is defined using Dirichlet series, while the Gothic fonts used for $ \mathfrak{S}^{r,G}_{\cusp}$ stands for taking residue at $s=1$.

\bigskip
Langlands' idea is to establish the existence of the limits $S^{r,G}_{\cusp}$ (resp. $\mathfrak{S}^{r,G}_{\cusp}$) directly from studying the geometric side (and concurrently the spectral side) of the stable trace formula, and to establish the beyond endoscopic decomposition of $S^{r,G}_{\cusp}$ (resp. $\mathfrak{S}^{r,G}_{\cusp}$) that reflects the order of poles at $s=1$ of Langlands $L$-functions $L(s,\pi,r)$ with respect to $r$, attached to discrete automorphic represenations $\pi$ of $G(\mathbf{A}_F)$ of Ramanujan type; this is a very difficult problem and we will refer to section 5 for some related discussions. In this section, we would like to verify that the limit $\mathfrak{S}^{r,G}_{\cusp}(f)$ exists for $r=r_1,r_2$, and also to establish $\star$-beyond endoscopic decomposition for $\mathfrak{S}^{r,G}_{\cusp}$ similar to the case of $S^G_{\cusp}$; for $r=r_1$ one in fact has the vanishing of $\mathfrak{S}^{r_1,G}_{\cusp}$. The final result is given in Theorem 4.1 below. In the case of $S^{r,G}_{\cusp}(f)$, we can establish the existence of the limit under Hypothesis 4.2 on the test function $f$ (Theorem 4.4 below).

\bigskip

We first make the definiton: for each elliptic $\star$-beyond endoscopic data $(H,\rho)$ of $G$, define, following [LP], the {\it dimension data} with respect to $(H,\rho)$ to be the following function on the set of algebraic representations $r$ of $\widehat{G}$:
\begin{eqnarray*}
m_{(H,\rho)}(r) : = \dim_{\mathbf{C}} \Hom_{\widehat{H}}(1,r \circ \rho )
\end{eqnarray*}
Clearly we have $m_{(H,\rho)}(r_1) =0$ for each elliptic $\star$-beyond endoscopic data $(H,\rho)$. On the other hand, for $r=r_2$, we have by Proposition 3.1 of [J]:
\begin{eqnarray}
m_{(H,\rho)}(r_2) = k-1
\end{eqnarray}
if $H$ is as in (3.1). In fact, by the results in section 4 of [J], for an elliptic $\star$-beyond endoscopic data $(H,\rho)$ for $G$ as in (3.1), the partition
\[
2N =\sum_{i=1}^k m_i
\]
can be recovered from the set of dimension data $m_{(H,\rho)}(r_a)$ for $a=2,4,\cdots,2[N/2]$ (for $a$ odd one has $m_{(H,\rho)}(r_a)=0$).

\bigskip
One would like to establish the existence of the limits $S^{r_a,G}_{\cusp},\mathfrak{S}^{r_a,G}_{\cusp}$ and also their $\star$-beyond endoscopic decomposition. More precisely, for each elliptic $\star$-beyond endoscopic data $(H,\rho)$, and each algebraic representation $r$ of $\widehat{G}$, and $n \geq 1$, we can define the function $\leftexp{n}{h}^{r \circ \rho}_w \in \mathcal{H}^{sph}(H(F_v))$ by a similar condition as in (4.1), namely that its Satake transform $\widehat{\leftexp{n}{h}^{r \circ \rho}_w}$ satisfies
\begin{eqnarray}
\widehat{ \leftexp{n}{h}^{r \circ \rho}_w}(d) = tr ((r \circ \rho)(d^n) )
\end{eqnarray}
for each semi-simple conjugacy class $d$ in $\widehat{H}$. Again we denote the function simply as $h_w^{r \circ \rho}$ when $n=1$. We then have, under the algebra homomorphism $b$ as in (3.11), the following compatibility:
\begin{eqnarray}
b(\leftexp{n}{h}^{r}_w) =\leftexp{n}{h}^{r \circ \rho}_w
\end{eqnarray}

\bigskip

Then for each test function $f^{\prime} \in \mathcal{H}(H(\mathbf{A}_F))$ spherical at places outside $S$, one would like to establish the existence of the limits, for $a=1,2$:
\begin{eqnarray}
& & \leftexp{\star}{P}^{r_a \circ \rho, H}_{\cusp}(f^{\prime}) \\
&:=& \lim_{X \rightarrow \infty} \frac{1}{X} \sum_{Nw \leq X, w \notin S} \log Nw \cdot \leftexp{\star}{P}^{H}_{\cusp}((f^{\prime})^{r_a \circ,\rho,w}) \nonumber
\end{eqnarray}
\begin{eqnarray}
 & & \leftexp{\star}{\mathfrak{P}}^{r_a \circ \rho, H}_{\cusp}(f^{\prime}) := \lim_{s=1} (s-1) \cdot \leftexp{\star}{\mathcal{P}}^{r_a \circ \rho,H}_{\cusp} (f^{\prime},s) \\
& & \leftexp{\star}{\mathcal{P}}^{r_a \circ \rho,H}_{\cusp}(f^{\prime},s) :=  \sum_{ w \notin S} \sum_{n \geq 1} \frac{\log Nw}{Nw^{ns}} \cdot \leftexp{\star}{P}^{H}_{\cusp}(\leftexp{n}{(f^{\prime})}^{r_a \circ,\rho,w})  \nonumber
\end{eqnarray}
with a similar explanation of the notations as before, and to show that
\begin{eqnarray}
& & \leftexp{\star}{P}^{r_a \circ \rho, H}_{\cusp}(f^{\prime}) = m_{(H,\rho)}(r_a) \cdot \leftexp{\star}{P}^{H}_{\cusp}(f^{\prime}) \\
& &  \leftexp{\star}{\mathfrak{P}}^{r_a \circ \rho, H}_{\cusp}(f^{\prime}) = m_{(H,\rho)}(r_a) \cdot \leftexp{\star}{P}^{H}_{\cusp}(f^{\prime}) \nonumber
\end{eqnarray}
from which both the existence of the limits, and also the $\star$-beyond endoscopic decompositions, for $S^{r_a,G}_{\cusp},\mathfrak{S}^{r_a,G}_{\cusp}$ ($a=1,2$), would follow.

\bigskip
We first consider the case for $\leftexp{\star}{\mathfrak{P}}^{r_a \circ \rho,H}_{\cusp}(f^{\prime})$. To do this, one would like to use the known results on the Langlands' $L$-function. Consider a parameter $\phi_H \in \Phi^{\tau}_{\prim}(H)$ as in (3.5) that occurs in the definition of $\leftexp{\star}{P}^H_{\cusp}(f^{\prime})$. Since the test function $f^{\prime}$ is spherical at places outside $S$, we may assume that $\phi_H$ belongs to the subset $\Phi^{S,\tau}_{\prim}(H)$ of $\Phi^{\tau}_{\prim}(H)$ consisting of parameters that is unramified outside $S$ (as otherwise its contribution to $\leftexp{\star}{P}^H_{\cusp}(f^{\prime})$ is zero), i.e. for each place $w \notin S$, the localization $\phi_{H,w}$ of $\phi_H$ at $w$ is unramified, and hence has the associated Frobenius-Hecke conjugacy class $c(\phi_{H,w})$ in $\widehat{H}$, with the local $L$-factor
\[
L(s,\phi_{H,w},r_a \circ \rho) =\frac{1}{\det(I - (r_a \circ \rho)(c(\phi_{H,w})))}
\]
and the partial $L$-function is initially defined as an Euler product:
\[
L^S(s,\phi_{H},r_a \circ \rho) =\prod_{w \notin S} L(s,\phi_{H,w},r_a \circ \rho) 
\]
which converges absolutely in some right half plane $Re(s) \gg 0$; using the $1/2$ weak Ramanujan bound of Jacquet-Shalika ({\it c.f.} Corollary 2.5 of [JS]), applied to the components $\phi_1,\cdots,\phi_k$ of $\phi_H$, we see that the Euler product for $L^S(s,\phi_{H},r_a \circ \rho)$ converges absolutely in the region $Re(s) > 1 + a/2$ (one can of course invoke slightly stronger bounds as given for example in Theorem 2 of [LRS], but this will not be needed in the sequel). In addition, this is also equal to the partial $L$-function 
\[
L^S(s,\phi,r_a)
\]
for the parameter $\phi := \rho \circ \phi_H \in \Phi^{\tau}_2(G)$, with respect to $r_a$. Here we recall that the Langlands $L$-function $L(s,\pi,r_a)$, for irreducible unitary representations $\pi$ of $G(\mathbf{A}_F)$ that belong to the same global Arthur packet classified by the parameter $\phi$, are all equal, hence the notation $L^S(s,\phi,r_a)$ is justified (similarly for $L^S(s,\phi_H,r_a \circ \rho)$). Finally, by the local Langlands correspondence for $G$ as established in chapter 6 of [A1], one can in fact defined the local $L$-factors at the places $w \in S$, and hence the completed $L$-function, but this is not necessary for our purpose.

\bigskip

In the rest of this section, we assume that $a=1,2$. Then the partial $L$-functions $L^S(s,\phi_H,r_a \circ \rho)=L^S(s,\phi,r_a)$ all have meromorphic continuation to the whole of $\mathbf{C}$, as we explain below.

\bigskip

\bigskip
Firstly, denote by $\Pi$ the strong functorial transfer of $\phi_H$ to $GL(2N,\mathbf{A}_F)$ with respect to the embedding $std \circ \rho$ of $\widehat{H}$ into $GL(2N,\mathbf{C})$, namely for the parameter $\phi_H$ as in (3.5), then $\Pi$ is the isobaric direct sum of $\phi_1,\cdots,\phi_k$. One then has:
\begin{eqnarray*}
L^S(s,\phi_H,r_1 \circ \rho) =L^S(s,\phi,r_1)=L(s,\Pi)=\prod_{i=1}^k L^S(s,\phi_i)
\end{eqnarray*}
As for $r_2$, we follow section 2 of [J] ({\it c.f.} Theorem 2.3 of {\it loc. cit.}): using the decomposition $\Lambda^2 std = r_2 \oplus 1$, one has
\begin{eqnarray*}
L^S(s,\phi_H,r_2,\circ \rho) =L^S(s,\phi,r_2)=L^S(s,\Pi, \Lambda^2)/\xi^S_F(s)
\end{eqnarray*}
where $L^S(s,\Pi, \Lambda^2)$ is the partial exterior square $L$-function of $\Pi$, and $\xi^S_F(s)$ is the partial Dedekind zeta function of $F$; in addition we have:
\begin{eqnarray*}
L^S(s,\Pi,\Lambda^2) = \prod_{i=1}^k L^S(s,\phi_i,\Lambda^2) \cdot \prod_{1\leq i<j\leq k } L^S(s,\phi_i \times \phi_j)
\end{eqnarray*}
with $L^S(s,\phi_i \times \phi_j)$ being the partial Rankin-Selberg $L$-functions for the pair $\phi_i,\phi_j$. 

\bigskip

All the $L$-functions that occur have meromorphic continuation to the whole of $\mathbf{C}$.

\bigskip

The analytic continuation properties of $\xi_F^S(s)$ is of course well known. It is meromorphic on $\mathbf{C}$ with only a simple pole at $s=1$, and the Euler product is absolutely convergent in the region $Re(s) > 1$, hence non-vanishing in this region.

\bigskip
The functions $L^S(s,\phi_i)$ are the principal $L$-functions of Godement--Jacquet [GJ], and are entire on the whole of $\mathbf{C}$, and by [JS1] are non-vanishing on the line $Re(s) = 1$; in addition the Euler product being absolutely convergent in the region $Re(s) >1$, by Theorem 5.3 of [JS], and hence non-vanishing in the region. Thus they are holomorphic and non-vanishing in the region $Re(s) \geq 1$; hence the same is true for $L^S(s,\Pi)=L^S(s,\phi_H,r_1 \circ \rho)$; in particular the order of pole of $L^S(s,\phi_H,r_1 \circ \rho)$ at $s=1$ is equal to zero, i.e. equal to $m_{(H,\rho)}(r_1)$. Here as a side remark, we note that from the proof of Theorem 5.3 of [JS], one has the stronger result that the Dirichlet series for $\log L^S(s,\phi_i)$:
\[
\sum_{w \notin S} \sum_{n \geq 1} \frac{1}{ n \cdot Nw^{ns}} \cdot tr(c(\phi_{i,w})^n)
\]
is absolutely convergent in the region $Re(s)>1$, from which we deduce that the Dirichlet series for the logarithmic derivative $-d/ds \log L^S(s,\phi_i)$
\[
\sum_{w \notin S} \sum_{n \geq 1} \frac{\log Nw}{ Nw^{ns}} \cdot tr(c(\phi_{i,w})^n)
\]
is also absolutely convergent in the region $Re(s) >1$, and hence that the Dirichlet series for the logarithmic derivative of $L^S(s,\phi_H,r_1 \circ \rho)$ is absolutely convergent in the region $Re(s)>1$ as well.

\bigskip
As for $L^S(s,\phi_H,r_2,\circ \rho)$, we use the known analytic properties of the partial $L$-functions $L^S(s,\phi_i,\Lambda^2), L^S(s,\phi_i \times \phi_j)$. Firstly, by Shahidi [Sh], they are meromorphic on the whole of $\mathbf{C}$ (for the $L^S(s,\phi_i \times \phi_j)$ one can also use [CPS]). In addition, again by Theorem 5.3 of [JS], the Euler product defining the Rankin-Selberg L-functions $L^S(s,\phi_i \times \phi_j)$ (here for the moment we allow $i=j$) are absolutely convergent in the region $Re(s) >1$ (again the proof in {\it loc. cit.} actually shows that the series for the logarithmic derivative of $L^S(s,\phi_i \times \phi_j)$ is absolutely convergent in the region $Re(s)>1$), and thus the Rankin-Selberg $L$-functions are non-vanishing in this region. 

\bigskip
Now we note in the particular case $i=j$, we have the factorization:
\[
L^S(s,\phi_i \times \phi_i ) =  L^S(s,\phi_i,\Lambda^2) \cdot L^S(s,\phi_i,\Sym^2)
\]
One knows, from the work of D. Belt [B], and Takeda [T1,T2] respectively, that the exterior square $L$-function $L^S(s,\phi_i,\Lambda^2)$ and the symmetric square $L$-function $L^S(s,\phi_i,\Sym^2)$ are holomorphic in the region $Re(s)>1$ (here we remark that both authors show that the functions $L^S(s,\phi_i,\Lambda^2)$ and respectively $L^S(s,\phi_i,\Sym^2)$ can have poles only at $s=0,1$, but we do not need this stronger fact). Thus since the Rankin-Selberg $L$-function $L^S(s,\phi_i \times \phi_i)$ are holomorphic and non-vanishing in this region, we have the same non-vanishing for $L^S(s,\phi_i,\Lambda^2)$ (and $L^S(s,\phi_i,\Sym^2)$) in the region $Re(s)>1$. Note that the non-vanishing of $L^S(s,\phi_i,\Lambda^2)$ (and $L^S(s,\phi_i,\Sym^2)$) in the region $Re(s) >1$ is established without knowing the absolute convergence of their Euler product in the region $Re(s)>1$ (indeed this appear to be not available in the literature, unless one has the exact Ramanujan bound for $\phi_i$).

\bigskip
On the line $Re(s)=1$, we know from Theorem 5.1 [Sh] that the $L$-functions $L^S(s,\phi_i \times \phi_j)$ (here $i \neq j$) and $L^S(s,\phi_i,\Lambda^2)$ (and $ L^S(s,\phi_i,\Sym^2)$) are holomorphic and non-vanishing, with the exception that the exterior square $L$-functions $L^S(s,\phi_i,\Lambda^2)$ have a simple pole at $s=1$ (recall that the $\phi_i$'s are of symplectic type), while the Rankin-Selberg $L$-functions $ L^S(s,\phi_i \times \phi_j)$ do not have poles at $s=1$, since for $i \neq j$ the $\phi_i,\phi_j$ are inequivalent (in fact the functions $L^S(s,\phi_i \times \phi_j)$ are entire on $\mathbf{C}$, by [CPS]).

\bigskip
Hence one conclude that $L^S(s,\phi_H,r_2 \circ \rho)$ is meromorphic on the whole of $\mathbf{C}$; it is holomorphic and non-vanishing in the region $Re(s) \geq 1$, {\it except} that it has a possible pole at $s=1$, of order equal to
\[
k-1 = m_{(H,\rho)}(r_2)
\]

Thus the logarithmic derivative 
\begin{eqnarray}
 - \frac{d}{ds} \log L^S(s,\phi_H,r_a \circ \rho), \,\ a=1,2
\end{eqnarray}
is meromorphic on the whole of $\mathbf{C}$; it is holomorphic in the region $Re(s) \geq 1$, except for a possible simple pole at $s=1$, with residue equal to the order of pole of $L^S(s,\phi_H,r_a\circ \rho) $ at $s=1$, in other words equal to $m_{(H,\rho)}(r_a)$.
\bigskip

Now for the logarithmic derivative $- d/ds \log L^S(s,\phi_H,r_a \circ \rho)$, its Dirichlet series is given as:
\begin{eqnarray}
  \sum_{w \notin S} \sum_{n=1}^{\infty} \frac{log Nw}{Nw^{ns}} \cdot tr ((r_a \circ \rho)(c(\phi_{H,w})^n) ) \end{eqnarray}
which converges absolutely in the region $Re(s)>1 + a/2$, where as before, the region of absolute convergence comes from the $1/2$ weak Ramanujan bound of [JS] (and also as noted earlier, converges absolutely in the region $Re(s) >1$ when $a=1$). 

\bigskip

We thus have, in the region $Re(s) > 1+ a/2$:
\begin{eqnarray*}
& & \leftexp{\star}{\mathcal{P}}^{r_a \circ \rho,H}_{\cusp}(f^{\prime},s) \\
&=& \sum_{w \notin S} \sum_{n \geq 1} \frac{\log Nw}{Nw^{ns}} \cdot \leftexp{\star}{P}^{H}_{\cusp}(  (\leftexp{n}{(f^{\prime})}^{r_a \circ \rho,w} ) \\
&=& \sum_{w \notin S} \sum_{n \geq 1} \frac{\log Nw}{Nw^{ns}}  \sum_{\phi_H \in \Phi_{\prim}(H)} 
(\leftexp{n}{(f^{\prime})}^{r_a \circ \rho,w})^H(\phi_H) \\
&=& \sum_{w \notin S} \sum_{n \geq 1} \frac{\log Nw}{Nw^{ns}}  \sum_{\phi_H \in \Phi_{\prim}(H)} \widehat{ \leftexp{n}{h}^{r_a \circ \rho,w}  }(c(\phi_{H,w}))   \cdot (f^{\prime})^H(\phi_H) \\
&=&   \sum_{\phi_H \in \Phi_{\prim}(H)}  \Big( \sum_{w \notin S} \sum_{n \geq 1} \frac{\log Nw}{Nw^{ns}} \cdot tr(( r_a \circ \rho)( c(  \phi_{H,w} )^n  )) \Big) \cdot (f^{\prime})^H(\phi_H) \\
&=&   \sum_{\phi_H \in \Phi_{\prim}(H)} - \frac{d}{ds} \log L^S(s,\phi_H,r_a \circ \rho) \cdot (f^{\prime})^H(\phi_H)
\end{eqnarray*}
Here the interchange of the sums is justified, since by the admissibility property of $\leftexp{\star}{P}^H_{\cusp}$, there is a finite subset $\Phi^{\prime}$ of $\Phi^{S,\tau}_{\prim}(H)$ (independent of $w \notin S$ and $n \geq 1$), such that for any $w \notin S$, the summands occurring in the definition of $\leftexp{\star}{P}^H_{\cusp}(\leftexp{n}{(f^{\prime})}^{r_a\circ \rho,w})$ vanishes for $\phi_H \notin \Phi^{\prime}$. Thus in the above, the sum over $\phi_H \in \Phi_{\prim}(H)$ can be taken over the finite set $\Phi^{\prime}$.

\bigskip
In particular $\leftexp{\star}{\mathcal{P}}^{r_a,H}_{\cusp}(f^{\prime},s) $ is absolutely convergent in the region $Re(s) > 1+ a/2$, and has holomorphic analytic continuation in the region $Re(s) \geq 1$ except for a possible pole at $s=1$. In addition:
\begin{eqnarray*}
& & \leftexp{\star}{\mathfrak{P}}^{r_a \circ \rho,H}_{\cusp}( f^{\prime} ) \\
&=& \lim_{s \rightarrow 1} (s-1) \cdot \mathcal{P}^{r_a \circ \rho,H}_{\cusp}(f^{\prime},s) \\
&=&  \sum_{\phi_H \in \Phi_{\prim}(H)}  \lim_{s \rightarrow 1 }-(s-1) \cdot  \frac{d}{ds} \log L^S(s,\phi_H,r_a \circ \rho) \cdot (f^{\prime})^H(\phi_H) \\
&=&   \sum_{\phi_H \in \Phi_{\prim}(H)}  m_{(H,\rho)}(r_a) \cdot (f^{\prime})^H(\phi_H) \\
&=& m_{(H,\rho)}(r_a)  \cdot \leftexp{\star}{P}^{H}_{\cusp}(f^{\prime})
\end{eqnarray*}

\bigskip
Finally back to $G$. Given test function $f \in \mathcal{H}(G(\mathbf{A}_F))$, choose, for each $(H,\rho)$, a test function $f^{\prime} \in \mathcal{H}(H(\mathbf{A}_F))$ that corresponds to $f$ under the stable transfer, thus
\[
f^H=(f^{\prime})^H
\] 
In addition recall that the stable transfer $f \rightarrow f^H$ is compatible with the algebra homomorphism $b$ of (3.11) at the spherical places, and we have the compatibility that $b(\leftexp{n}{h}^{r_a}_w) =\leftexp{n}{h}^{r_a \circ \rho}_w$. Thus we have:
\[
(\leftexp{n}{f}^{r_a,w})^H =( \leftexp{n}{(f^{\prime})}^{r_a \circ \rho,w})^H
\]
It follows, on applying Theorem 3.6 to $S^G_{\cusp}(\leftexp{n}{f}^{r_a,w})$ for each $n \geq 1$ and $w \notin S$, that we have:
\begin{eqnarray*}
& & \mathcal{S}^{r_a,G}_{\cusp}(f) = \sum_{w \notin S} \sum_{n \geq 1} \frac{\log Nw}{Nw^s} \cdot S^G_{\cusp}(\leftexp{n}{f}^{r_a,w}   ) \\
&=&  \sum_{w \notin S} \sum_{n \geq 1} \frac{\log Nw}{Nw^s} \sum_{(H,\rho)}  \iota(G,H) \widehat{\leftexp{\star}{P}^H_{\cusp}}( (\leftexp{n}{f}^{r_a,w})^H   ) \\
&=& \sum_{(H,\rho)} \iota(G,H)   \sum_{w \notin S} \sum_{n \geq 1} \frac{\log Nw}{Nw^s}\cdot  \leftexp{\star}{P}^H_{\cusp}( \leftexp{n}{(f^{\prime})}^{r_a \circ \rho,w})  \\
&=&  \sum_{(H,\rho)} \iota(G,H)   \cdot \leftexp{\star}{\mathcal{P}}^{r_a \circ \rho,H}_{\cusp}(f^{\prime},s)
\end{eqnarray*}
is absolutely convergent in the region $Re(s)>1+ a/2$, has holomorphic analytic continuation in the region $Re(s) \geq1$, except for a possible pole at $s=1$. Upon taking residue at $s=1$, we see that the limit $\mathfrak{S}^{r_a,G}_{\cusp}(f)$ exists, and we have:
\begin{eqnarray*}
\mathfrak{S}^{r_a,G}_{\cusp}(f) = \sum_{(H,\rho)} \iota(G,H) \cdot m_{(H,\rho)}(r_a) \cdot \leftexp{\star}{P}^{H}_{\cusp}(f^{\prime})
\end{eqnarray*}

\bigskip

Define $\iota(r_a,H):= \iota(G,H) \cdot m_{(H,\rho)}(r_a)$. Then to sum up, we have:
\begin{theorem}
For an elliptic $\star$-beyond endoscopic data $(H,\rho)$ of $G$, test function $f^{\prime} \in \mathcal{H}(H(\mathbf{A}_F))$, and $a=1,2$, the Dirichlet series $\leftexp{\star}{\mathcal{P}}^{r_a,H}_{\cusp}(s,f^{\prime})$ converges absolutely in the region $Re(s)>1 + a/2$, and has holomorphic analytic continuation in the region $Re(s)>1$. The limit $\leftexp{\star}{\mathfrak{P}}^{r_a,H}_{\cusp}(f^{\prime})$ exists, and one has
\[
\leftexp{\star}{\mathfrak{P}}^{r_a,H}_{\cusp}(f^{\prime}) = m_{(H,\rho)}(r_a) \cdot \leftexp{\star}{P}^H_{\cusp}(f^{\prime})
\]
Consequently, for any $f \in \mathcal{H}(G(\mathbf{A}_F))$ the Dirichlet series $\mathcal{S}^{r_a,G}_{\cusp}(f)$ converges absolutely for $Re(s) > 1+ a/2$, has holomorphic analytic continuation in the region $Re(s)>1$. The limit $\mathfrak{S}^{r_a,G}_{\cusp}(f)$ exists, and we have the $\star$-beyond endoscopic decomposition:
\[
\mathfrak{S}^{r_a,G}_{\cusp}(f) = \sum_{(H,\rho)} \iota(r_a,H) \cdot \widehat{\leftexp{\star}{P}^{H}_{\cusp}}(f^H)
\]
In particular $\mathfrak{S}^{r_1,G}_{\cusp}(f)=0$.
\end{theorem}

\bigskip
Next we turn to $S^{r_a,G}_{\cusp}$. Following [A2] and also the afterword section of [A3], we would like to apply the Wiener-Ikehara Tauberian argument, as in section A.3 of chapter one in [Se]. However, the argument in {\it loc. cit.} requires one to know the exact Ramanujan bound. In Langlands' beyond endoscopy proposal, one hopes to establish the Generalized Ramanujan Conjecture as a consequence of the Principle of Functoriality. In this paper, our purpose is more modest and only aims to give examples for which the $r$-stable trace formula is valid. So we will make use of known results on the Generalized Ramanujan Conjecture. Hence we make the following

\begin{hypothesis}
The field $F$ is totally real, and for each archimedean place $v$ of $F$, the component $f_v$ of the test function $f \in \mathcal{H}(G(\mathbf{A}_F))$ at $v$, is a linear combination of stable pseudo-coefficients for $L$-packets of discrete series representations, whose infinitesimal characters are strongly regular.
\end{hypothesis}

We now explain the condition. When $F$ is totally real, we have, for each archimedean place $v$, the bijection between the set of $L$-packets of discrete series representations of $G(F_v) (\cong G(\mathbf{R}))$, with the set of irreducible finite dimensional representations of $G(F_v)$ with unitary central character, as in [A8]. Thus given $\mu$, we denote by $\Pi(\mu)$ the corresponding $L$-packet of discrete series representations of $G(F_v)$, which consists of discrete series representation of $G(F_v)$ having the same infinitesimal character and central character as $\mu$. Also denote by $f_{v}(\mu)$ the stable pseudo-coefficient of the packet $\Pi(\mu)$. We then say that the infinitesimal character of the packet $\Pi(\mu)$ is strongly regular, if the highest weight of $\mu$ is non-singular. 

\bigskip
By linearity we assume without loss of generality that the test function $f$ is chosen so that $f_v$ is of the form $f_v(\mu_v)$ for some $\mu_v$ whose highest weight is non-singular. Then for any $\phi \in \Phi^{\tau}_2(G)$, if
\[
f^G(\phi) = \prod_v f_v^G(\phi_v) \neq 0
\] 
then in particular $f_v^G(\phi_v) \neq 0$ for each archimedean place $v$. By the argument of p. 283 of [A8] (relying on the result of Vogan-Zuckerman [VZ]), the non-singularity condition on the highest weight $\mu_v$ then allows us to deduce that
\begin{property} 
For each archimedean place $v$, the local parameter $\phi_v$ corresponds to the packet $\Pi(\mu_v)$. 
\end{property}
Thus only those $\phi$ satisfying Property 4.3 would contribute to $S^G_{\cusp}(f)$. Consequently, for each $(H,\rho)$, and $f^{\prime} \in \mathcal{H}(H(\mathbf{A}_F))$ corresponding to $f$ under the stable transfer, only those $\phi_H \in \Phi_{\prim}(H)$ such that $\phi = \rho \circ \phi_H$ satisfies Property 4.3, would contribute to $P^H_{\cusp}(f^{\prime})$, and so we only need to deal with those $\phi$ (and $\phi_H$) that satisfies Property 4.3. Write the parameter $\phi$ in the form
\[
\phi = \phi_1 \boxplus \cdots \boxplus \phi_k
\]
as in (2.1). Each $\phi_i$ is then a self-dual cuspidal automorphic representation of $GL(m_i,\mathbf{A}_F$) (of symplectic type); and each archimedean local component $\phi_{i,v}$ is cohomological. As $F$ is totally real, the exact Ramanujan bound for $\phi_i$ is known to hold; in our present case the archimedean local component $\phi_{i,v}$ satisfies the regularity condition of [Shi], and the results of [Shi] actually show that $\phi_i$ is tempered at every place of $F$. Thus the same is true for $\phi$ (here we remark that the main results of [Shi] apply in the setting of conjugate self-dual cuspidal automorphic representations of $GL(m,\mathbf{A}_K)$ of cohomological type, for $K$ a CM field, but the Ramanujan Conjecture for self-dual cuspidal automorphic representations on $GL(m,\mathbf{A}_F)$ of cohomological type, for totally real $F$ can be reduced to the conjugate self-dual case by a standard base change argument, as given for example in section 1 of [C]).

\bigskip
With the exact Ramanujan bound, one then has the much stronger inequality:
\[
| tr (r_a \circ \rho)(c(\phi_{H,w})^n) |  \leq \dim r_a 
\]

We decompose (4.11) as the sum
\[
F_a(s,\phi_H) +E_a(s,\phi_H) 
\]
where
\begin{eqnarray*}
F_a(s,\phi_H)&=&  \sum_{w \notin S} \frac{log Nw }{ Nw^s}  \cdot tr( (r_a \circ \rho)(c(\phi_{H,w})))\\
E_a(s,\phi_H) &=& \sum_{w \notin S} \sum_{n=2}^{\infty} \frac{log Nw }{ Nw^{ns}}  \cdot tr( (r_a \circ \rho)(c(\phi_{H,w})^n))
\end{eqnarray*}

We follow the argument of section A.3 of chapter one of [Se]. The Dirichlet series for $F_a(s,\phi_H)$ converges absolutely for $Re(s)>1$, and whose coefficients are dominated by that of the series
\[
\dim r_a \cdot \sum_{w \notin S} \frac{\log Nw}{ Nw^s}  
\]
while $E_a(s,\phi_H)$ converges absolutely in the region $Re(s) >1/2$. Thus the analytic behaviour of $F_a(s,\phi_H)$ in the region $Re(s) \geq 1$ is the same as that of the logarithmic derivative (4.10), which we have already established. Thus $F_a(s,\phi_H)$ satisfies the hypothesis of the Wiener-Ikehara Tauberian theorem ({\it c.f. loc. cit.}), and so we deduce that:

\begin{eqnarray*}
 & &  \lim_{X \rightarrow \infty} \frac{1}{X} \sum_{Nw \leq X,w \notin S }  \log Nw \cdot \tr( (r_a \circ \rho)(c(\phi_{H,w})))\\  &=& \Res_{s=1} F_a(s,\phi_H) =   \Res_{s=1} \Big( -\frac{d}{ds}  \log L^S(s,\phi_H,r_a \circ \rho)  \Big)  \\
&=& m_{(H,\rho)}(r_a)
\end{eqnarray*}
Thus the limit defining $\leftexp{\star}{P}^{r_a \circ \rho,H}_{\cusp}(f^{\prime}) $ is equal to:
\begin{eqnarray*}
& & \lim_{X \rightarrow \infty} \frac{1}{X} \sum_{Nw \leq X, w \notin S} \log Nw \cdot \leftexp{\star}{P}^{H}_{\cusp}((f^{\prime})^{r_a \circ \rho,w})\\
&=&  \lim_{X \rightarrow \infty} \frac{1}{X} \sum_{Nw \leq X, w \notin S} \sum_{\phi_H \in \Phi^{\tau}_{\prim}(H)} \log Nw \cdot (f^{\prime})^{r_a \circ \rho,w}(\phi_H) \\
&=&  \lim_{X \rightarrow \infty} \frac{1}{X} \sum_{Nw \leq X, w \notin S} \sum_{\phi_H \in \Phi^{\tau}_{\prim}(H)}  \log Nw \cdot \widehat{h}^{r_a \circ \rho}(c(\phi_{H,w})) \cdot (f^{\prime})^H(\phi_H) \\
&=&  \sum_{\phi_H \in \Phi^{\tau}_{\prim}(H)} \Big( \lim_{X \rightarrow \infty} \frac{1}{X} \sum_{Nw \leq X, w \notin S}  \log Nw \cdot tr((r_a \circ \rho)(c(\phi_{H,w}))) \Big) \cdot (f^{\prime})^H(\phi_H) \\
&=& m_{(H,\rho)}(r_a) \cdot \leftexp{\star}{P}_{\cusp}^H(f^{\prime})
\end{eqnarray*}
with the interchange of the limit and sum being justified as in the case for $\leftexp{\star}{\mathfrak{P}}^{r_a,G}_{\cusp}$.  

\bigskip

Finally back to $G$. Again for each $(H,\rho)$, if we pick $f^{\prime} \in \mathcal{H}(H(\mathbf{A}_F))$ such that
\[
f^H=(f^{\prime})^H
\] 
so that
\[
(f^{r_a,w})^H =( (f^{\prime})^{r_a \circ \rho,w})^H
\]
then we have
\begin{eqnarray*}
& & \lim_{X \rightarrow \infty }\frac{1}{X} \sum_{Nw \leq X, w \notin S}  \log Nw \cdot S^{G}_{\cusp}(f^{r_a,w}) \\
&=&  \lim_{X \rightarrow \infty }\frac{1}{X} \sum_{Nw \leq X, w \notin S} \log Nw \sum_{(H,\rho)} \iota(G,H) \cdot \widehat{\leftexp{\star}{P}^{H}_{\cusp}}((f^{r_a,w})^H)\\
&=& \sum_{(H,\rho)} \iota(G,H) \lim_{X \rightarrow \infty }\frac{1}{X} \sum_{Nw \leq X, w \notin S} \log Nw \cdot  \leftexp{\star}{P}^{H}_{\cusp}((f^{\prime})^{r_a\circ \rho,w}))\\
&=& \sum_{(H,\rho)} \iota(r_a,H) \cdot \leftexp{\star}{P}^H_{\cusp}(f^{\prime})
 \end{eqnarray*}
as required.

\bigskip
To sum up we have
\begin{theorem}
Assume Hypothesis 4.2 on the test function $f \in \mathcal{H}(G(\mathbf{A}_F))$ (in particular $F$ is totally real). Then the limit $S^{r_a,G}_{\cusp}(f)$, and we have the $\star$-beyond endoscopic decomposition, for $a=1,2$:
\begin{eqnarray*}
S^{r_a,G}_{\cusp}(f) = \sum_{(H,\rho)} \iota(r_a,H) \cdot \widehat{\leftexp{\star}{P}^{H}_{\cusp}}(f^H)
\end{eqnarray*}
In particular $S^{r_1,G}_{\cusp}(f)=0$.
\end{theorem}

\section{Some questions}

In this final section, we would like to pose some questions related to the results of the paper.

\bigskip

{\bf I}. In this paper we work with the split special odd orthogonal groups. Nevertheless using the results in [A1] and [M], it should be possible to establish similar results for other families of quasi-split classical groups.

\bigskip

{\bf II}. In this paper, we establish the vanishing of $\mathfrak{S}^{r_1,G}_{\cusp}(f)$, by combining results in [A1] and known results on the Langlands $L$-functions involved. In accordance with Langlands' beyond endoscopy ideas, one should try to establish this vanishing, but working directly with the geometric side of the stable trace formula. In the case of $G=SO(3)=PGL(2)$ over $\mathbf{Q}$, this is a consequence of the thesis of A. Altug [Alt] (under certain hypothesis on the test function $f$); in the process he also obtained the $1/4$-bound on the Ramanujan conjecture, purely from the method of Arthur-Selberg trace formula for $PGL(2)$. It would be very interesting to try to extend the technique of [Alt] and to establish the vanishing of $\mathfrak{S}^{r_1,G}_{\cusp}(f)$ for the rank two group $G=SO(5)=PGSp(4)$. In particular, one hopes that by working with the geometric side of the stable trace formula, one can establish the existence of the limits $\mathfrak{S}^{r_1,G}_{\cusp}(f), \mathfrak{S}^{r_2,G}_{\cusp}(f)$ without relying on the known results on the Langlands $L$-functions involved.

\bigskip
{\bf III}. In the case of $PGL(2)$, the only discrete automorphic representation that is not of Ramanujan type is the one-dimensional trivial representation, and in [Alt] this is removed from the trace formula by non-trivial application of the Poisson summation formula, following the suggestion of Frenkel-Langlands-Ng\^o [FLN]. For the group $PGSp(4)$, there are more discrete automorphic representations that are not of Ramanujan type. For a list of discrete automorphic representations of $PGSp(4)$ that are not of Ramanujan type, see [A7]. At present it is not clear how to extend the method of [FLN] to remove these other discrete automorphic representations that are not of Ramanujan type. A possible simplification is to consider test functions as in Hypothesis 4.2, namely whose archimedean components are stable pseudo-characters of $L$-packets of discrete series representations whose infinitesimal characters are sufficiently regular. Such a test function would then kill the contributions to the spectral side of the stable trace formula, from the discrete automorphic representations whose archimedean components are non-tempered (and as seen in section 4, the remaining discrete automorphic representations that contribute to the stable trace formula would then be of Ramanujan type). Also in this setting one has an explicit description of the geometric side of the stable trace formula [P] (which improves on the result [A8]).  

\bigskip

{\bf IV}. Theorem 4.1 applied to the case $G=SO(5)=PGSp(4)$, shows that the $r_2$-stable trace formula $\mathfrak{S}^{r_2,G}_{\cusp}$ precisely captures the discrete automorphic representations on $G$ that arise from $H=SO(3) \times SO(3) = PGL(2) \times PGL(2)$ (these are usually known as Yoshida lifts). Establishing the $\star$-beyond endoscopic decomposition for $\mathfrak{S}^{r_2,G}_{\cusp}$ in the case $G=SO(5)$, by directly working with the geometric side of the stable trace formula, would thus be an important case to test Langlands' beyond endoscopy ideas.

\bigskip

{\bf V}. We have see in section 4 that the partial $L$-function $L^S(s,\phi_H,r_2 \circ \rho) =L^S(s,\phi,r_2)$ for $\phi_H \in \Phi^{\tau}_{\prim}(H)$ is meromorphic on the whole of $\mathbf{C}$, holomorphic in the region $Re(s) \geq 1$, except possibly for a pole at $s=1$. Is possible to establish that the only poles of $L^S(s,\phi,r_2)$ are located at $s=0,1$? Such questions generally require the search for an integral representation of the $L$-function involved. In the case $N=2$, i.e. $G=SO(5)$, one can use the existence of functorial lifts with respect to $r_2$, of discrete automorphic representations of $SO(5)$ parametrized by $\phi$, to isobaric automorphic representations of $GL(5,\mathbf{A}_F$), using the argument of Theorem 6.6 of [K1] (this is a direct consequence of the exterior square lift from $GL(4)$ to $GL(6)$ of [K2], together with the tensor product lift of [R]). Thus for $G=SO(5)$, the only poles of $L^S(s,\phi,r_2)$ are at $s=0,1$, and the pole at $s=1$ occurs exactly when $\phi$ is of Yoshida type, i.e. $\phi=\phi_1 \boxplus \phi_2$, in which case the pole is simple.

\bigskip

{\bf VI}. The stable linear form $\leftexp{\star}{P}^G_{\cusp}$ defined as in (2.3), is of course only a first approximation to the primitive trace formula $P^G_{\cusp}$ as discussed in section 2 of [A2]. A possible inductive procedure to construct $P^G_{\cusp}$ is given for example in equation (2.6) of [A2]. From the point of view of global Arthur parameters, one would like to construct the subset $\Phi^{\tau}_{\pri}(G) \subset \Phi^{\tau}_{\simp}(G)$ of primitive parameters, namely those parameters that does not factor through any smaller subgroup of $\widehat{G}$. It is expected that the representations $\pi$ parametrized by the global Arthur packet attached to such a $\phi$, satisfies the following condition: for any irreducible algebraic representation $r$ of $\widehat{G}$ not equal to the trivial representation, the Langlands $L$-function $L(s,\pi,r)$ is holomorphic and non-vanishing in the region $Re(s)  \geq 1$, in particular has no poles at $s=1$; and conversely that given a discrete automorphic representation $\pi$ classified by a generic parameter $\phi \in \Phi^{\tau}_2(G)$, the condition on the $L$-function $L(s,\pi,r)$ for all such $r$  implies that $\phi \in \Phi^{\tau}_{\pri}(G)$. Then similar to (2.3), one would have:
\[
P^G_{\cusp}(f) =\sum_{\phi \in \Phi^{\tau}_{\pri}(G)} f^G(\phi)
\]
In addition, one expects that for any algebraic representation $r$ of $\widehat{G}$ that does not contain the trivial representation, the linear forms $P^{r,G}_{\cusp}, \mathfrak{P}^{r,G}_{\cusp}$, defined by the familiar limits:
\begin{eqnarray*}
& & P^{r,G}_{\cusp}(f) := \lim_{X \rightarrow \infty} \frac{1}{X} \sum_{Nw \leq X, w \notin S} \log Nw \cdot  P^{G}_{\cusp}(f^{r,w}) \\
& & \mathfrak{P}^{r,G}_{\cusp}(f) := \lim_{s \rightarrow 1} (s-1)\cdot \Big( \sum_{w \notin S} \sum_{n \geq 1} \frac{\log Nw}{Nw^{ns}} \cdot  P^{G}_{\cusp}(\leftexp{n}{f}^{r,w})  \Big)
\end{eqnarray*}
exist, and vanishes identically. 
\bigskip

{\bf VII}. The $L$-functions for $G$ with respect to the higher fundamental representations $r_a$, for $a >2$, in general are no longer related to $L$-functions of Langlands-Shahidi type, and there are also no known integral representations for the Langlands $L$-functions with respect $r_a$, so it is reasonable to expect that the study of the $r_a$-stable trace formula for $a>2$ would be much more involved. Nevertheless, the case $a=3$, i.e. the third fundamental representation $r_3$, and $G=SO(7), \widehat{G}=Sp(6,\mathbf{C})$, the Langlands $L$-function for $SO(7)$ respect to $r_3$ is of Langlands-Shahidi type (using the simply connected group of type $F_4$, {\it c.f.} [KS], Theorem 7.2.5). In addition, the poles of $L$-function for $SO(7)$ with respect to $r_3$ in the region $Re(s) \geq 1$, can be given in terms of the poles of exterior cube $L$-function for isobaric automorphic representations on $GL(6,\mathbf{A}_F)$ (since $\Lambda^3 std = r_3 \oplus std$). One also has the works of Ginzburg-Rallis [GR] and S. Yamana [Y] which give a complete description of poles of exterior cube $L$-functions for cuspidal automorphic representations on $GL(6,\mathbf{A}_F)$. 

\bigskip
More specifically, based on the results in [GR] and [Y], one would expect that the partial $L$-function $L^S(s,\phi,r_3)$ is entire, for $\phi$ a cuspidal parameter of $SO(7)$ (this is proved in [KS], Theorem 7.2.5, under suitable hypothesis). Consequently, one also expects that the $r_3$-stable trace formula $\mathfrak{S}^{r_3,G}_{\cusp}$ vanishes identically for $G=SO(7)$. It is certainly a very challenging problem to try to establish this from the geometric side of the stable trace formula.

\section{Appendix: Another construction of the stable transfer}

In this appendix, we give an alternative construction of the stable transfer as described in Theorem 3.6, which may be of independent interest.  

\bigskip

We need to use results on twisted endoscopy for which we refer to section 2,1 of [A1]. 

\bigskip

Firstly, from Corollary 2.1.2 of [A1], the Kottwitz-Shelstad transfer map 
\begin{eqnarray*}
\widetilde{\mathcal{H}}(2N,F_v) & \rightarrow & \mathcal{S}(G(F_v)) \\
\widetilde{f}_v & \mapsto & \widetilde{f}_v^G
\end{eqnarray*}
is surjective. Thus given $f_v \in \mathcal{H}(G(F_v))$, let $\widetilde{f}_v \in \widetilde{\mathcal{H}}(2N,F_v)$ be such that:
\[
\widetilde{f}_v^G = f_v^G. 
\]

Next we have the Levi embedding:
\[
GL(m_1) \times \cdots GL(m_k) \hookrightarrow GL(2N)
\]
such that the twisting automorphism of $GL(2N)$ preserves each of the factors $GL(m_i)$, and induces the corresponding twisting automorphism on the $GL(m_i)$ factors. Corresponding to this Levi embedding, we have the twisted descent map:
\begin{eqnarray*}
\widetilde{\mathcal{H}}(2N,F_v) \rightarrow \widetilde{\mathcal{H}}(m_1,F_v) \times \cdots \times \widetilde{\mathcal{H}}(m_k,F_v).
\end{eqnarray*}

Denote the image of $\widetilde{f}_v$ under the twisted descent map as $\widetilde{f}_{v,M}$. Finally we have the product of the Kottwitz-Shelstad transfer maps:
\begin{eqnarray}
& &  \widetilde{\mathcal{H}}(m_1,F_v) \times \cdots \times \widetilde{\mathcal{H}}(m_k,F_v) \rightarrow \\
& & \mathcal{S}(H_1(F_v)) \times \cdots \times  \mathcal{S}(H_k(F_v)) = \mathcal{S}(H(F_v)). \nonumber
\end{eqnarray}

Denote the image of $\widetilde{f}_{v,M}$ under (3.10) as $(\widetilde{f}_{v,M})^H$. We claim that this is the stable transfer of $f_v$ that is sought after.

\bigskip

Thus let $\phi_{H,v} = \phi_{1,v} \times \cdots \times \phi_{k,v}$ be a bounded local Langlands parameter of $H$ over $F_v$. For $i=1,\cdots,k$, denote by $\pi_{i,v}$ the irreducible tempered representation of $GL(m_i,F_v)$ corresponding to the parameter $\phi_{i,v}$ (regarded as a homomorphism $\phi_{i,v}: L_{F_v} \rightarrow GL(m_i,\mathbf{C})$). Similarly denote by $\Pi_v$ the irreducible tempered representation of $GL(2N,F_v)$ corresponding to the parameter $\rho \circ \phi_{H,v}$ (regarded as a homomorphism $\rho \circ \phi_{H,v}: L_{F_v} \rightarrow GL(2N,\mathbf{C})$); thus we may write $\rho \circ \phi_{H,v} = \phi_{1,v} \oplus \cdots \oplus \phi_{k,v}$). Then by the twisted character identities satisfied by tempered $L$-packets, stated for example as Theorem 2.2.1 of [A1], equation (2.2.3), we have:
\[
f_v^G(\rho \circ \phi_{H,v}) = \widetilde{f}_v^G(\rho \circ \phi_{H,v}) = \tr \widetilde{\Pi}_v(\widetilde{f}_v)
\]
where $\tr \widetilde{\Pi}_v$ is the twisted character of $\Pi_v$ with respect to Whittaker normalization ({\it c.f.} the discussion in section 2.2 of [A1], before the statement of Theorem 2.2.1 of {\it loc. cit.}) We will also use similar notation below for the twisted character $\tr \big( \widetilde{\pi}_{1,v} \times  \cdots \widetilde{\pi}_{k,v} \big) $ of $\pi_{1,v} \times \cdots \times \pi_{k,v}$.

\bigskip
On the other hand, the twisted descent map gives the identity:
\[
 \tr \widetilde{\Pi}_v(\widetilde{f}_v) = \tr  \big( \widetilde{\pi}_{1,v} \times  \cdots \widetilde{\pi}_{k,v} \big) ( \widetilde{f}_{v,M}).
\]

\bigskip
And finally, we apply again the twisted character identities given by Theorem 2.2.1 of [A1], equation (2.2.3), and obtain:
\[
\tr  \big( \widetilde{\pi}_{1,v} \times  \cdots \widetilde{\pi}_{k,v} \big) ( \widetilde{f}_{v,M}) = (\widetilde{f}_{v,M})^H(\phi_{1,v} \times \cdots \times \phi_{k,v}) =  (\widetilde{f}_{v,M})^H(\phi_{H,v}).
\]

\bigskip

Thus to conclude, given $f_v \in \mathcal{H}(G(F_v))$, by taking $f_v^H := (\widetilde{f}_{v,M})^H  \in \mathcal{S}(H(F_v))$, we have:
\[
f_v^H(\phi_{H,v}) = f_v^G(\rho \circ \phi_{H,v})
\]
for every bounded local Langlands parameters of $H$ over $F_v$. The map $f_v \mapsto f_v^H $ is thus the stable transfer as described in Theorem 3.6.

\end{document}